\documentclass[10pt,a4paper]{article}
\usepackage{amsfonts,amsmath}

\newtheorem{theorem}{Theorem}

\newtheorem{remark}{Remark}

\input{tcilatex}
\begin{document}

\date{}
\title{\textbf{A remark on the existence of solutions to a (}$k_{1}$,$k_{2}$%
\textbf{)-Hessian system with convection term}}
\author{ \vspace{1mm} \textsf{{Dragos-Patru Covei} } \\
{\small \textit{Department of Applied Mathematics}}\\
{\small \textit{Bucharest University of Economic Studies }}\\
{\small \textit{Piata Romana, 1st district, postal code: 010374, postal
office: 22, Romania}}\\
{\small \textit{e-mail address:coveidragos@yahoo.com}}}
\maketitle

\begin{abstract}
In this paper, two new results and short proofs are given for the existence
of positive entire large and bounded radial positive solutions for the
following nonlinear system with gradient term%
\begin{equation*}
\left\{ 
\begin{array}{l}
S_{k_{1}}\left( \lambda \left( D^{2}u_{1}\right) \right) +b_{1}\left(
\left\vert x\right\vert \right) \left\vert \nabla u_{1}\right\vert
^{k_{1}}=p_{1}\left( \left\vert x\right\vert \right) f_{1}\left(
u_{1},u_{2}\right) \text{ for }x\in \mathbb{R}^{N}\text{ , } \\ 
S_{k_{2}}\left( \lambda \left( D^{2}u_{2}\right) \right) +b_{2}\left(
\left\vert x\right\vert \right) \left\vert \nabla u_{2}\right\vert
^{k_{2}}=p_{2}\left( \left\vert x\right\vert \right) f_{2}\left(
u_{1},u_{2}\right) \text{ for }x\in \mathbb{R}^{N}\text{ ,}%
\end{array}%
\right. 
\end{equation*}%
where $S_{k_{i}}\left( \lambda \left( D^{2}u_{i}\right) \right) $ is the $%
k_{i}$-Hessian operator, $b_{1}$, $p_{1}$, $f_{1}$, $b_{2}$, $p_{2}$ and $%
f_{2}$ are continuous functions satisfying certain properties. Our results
complete and improve a recently work published by Zhang and Zhou. The main
difficulty in dealing with our system is the presence of the convection term.
\end{abstract}

\vspace{0.5cm}{\scriptsize 2010 AMS Subject Classification: Primary:35J25,
35J47 Secondary: 35J96.}


\section{Introduction}

The purpose of this paper is to present new results concerning the nonlinear
Hessian system with convection term 
\begin{equation}
\left\{ 
\begin{array}{l}
S_{k_{1}}\left( \lambda \left( D^{2}u_{1}\right) \right) +b_{1}\left(
\left\vert x\right\vert \right) \left\vert \nabla u_{1}\right\vert
^{k_{1}}=p_{1}\left( \left\vert x\right\vert \right) f_{1}\left(
u_{1},u_{2}\right) \text{ for }x\in \mathbb{R}^{N}\text{ (}N\geq 3\text{), }
\\ 
S_{k_{2}}\left( \lambda \left( D^{2}u_{2}\right) \right) +b_{2}\left(
\left\vert x\right\vert \right) \left\vert \nabla u_{2}\right\vert
^{k_{2}}=p_{2}\left( \left\vert x\right\vert \right) f_{2}\left(
u_{1},u_{2}\right) \text{ for }x\in \mathbb{R}^{N}\text{ (}N\geq 3\text{),}%
\end{array}%
\right.  \label{11}
\end{equation}%
where $b_{1}$, $p_{1}$, $f_{1}$, $b_{2}$, $p_{2}$, $f_{2}$ are continuous
functions satisfying certain properties, $k_{1},k_{2}\in \left\{
1,2,...,N\right\} $ and $S_{k_{i}}\left( \lambda \left( D^{2}u_{i}\right)
\right) $ stands for the $k_{i}$-Hessian operator defined as the sum of all $%
k_{i}\times k_{i}$ principal minors of the Hessian matrix $D^{2}u_{i}$. For
instance, the following well known operators are included in this class:%
\begin{equation*}
\begin{tabular}{llll}
\hline
\textbf{Operator:} & Laplacian &  & Monge--Amp\`{e}re \\ \hline
& $S_{1}\left( \lambda \left( D^{2}u_{i}\right) \right) =\Delta u_{i}=\func{%
div}\left( \nabla u_{i}\right) $ &  & $S_{N}\left( \lambda \left(
D^{2}u_{i}\right) \right) =\det \left( D^{2}u_{i}\right) $ \\ \hline
\end{tabular}%
\end{equation*}

In recent years equations of the type (\ref{11}) have been the subject of
rather deep investigations since appears from many branches of mathematics
and applied mathematics. For more surveys on these questions we advise the
paper of Alves and Holanda \cite{CA}, Bao-Ji and Li \cite{BAOII}, Bandle and
Giarrusso \cite{BA}, Cirstea and R\u{a}dulescu \cite{CR}, Cl\'{e}ment-Man%
\'{a}sevich and Mitidieri \cite{CL}, De Figueiredo and Jianfu \cite{FI},
Galaktionov and V\'{a}zquez \cite{GA},\ Ghergu and R\u{a}dulescu \cite{GR},
Salani \cite{SA}, Ji and Bao \cite{BAO}, Jian \cite{J}, Peterson and Wood 
\cite{PW}, Quittner \cite{Q}, Li and Yang \cite{LY}, Li-Zhang and Zhang \cite%
{Z}, Viaclovsky \cite{V,VI} and not the last Zhang and Zhou \cite{ZZ}.

The motivation for studying (\ref{11}) comes from the work of Ghergu and R%
\u{a}dulescu \cite{GR} where in discussion is a system of the type 
\begin{equation*}
\left\{ 
\begin{array}{l}
\Delta u_{1}+\left\vert \nabla u_{1}\right\vert =p_{1}\left( \left\vert
x\right\vert \right) f_{1}\left( u_{1},u_{2}\right) \text{ for }x\in \mathbb{%
R}^{N}\text{ (}N\geq 3\text{), } \\ 
\Delta u_{2}+\left\vert \nabla u_{2}\right\vert =p_{2}\left( \left\vert
x\right\vert \right) f_{2}\left( u_{1},u_{2}\right) \text{ for }x\in \mathbb{%
R}^{N}\text{ (}N\geq 3\text{),}%
\end{array}%
\right.
\end{equation*}%
and from the recently work of Zhang and Zhou \cite{ZZ} where the authors
have been considered the system%
\begin{equation*}
\left\{ 
\begin{array}{l}
S_{k}\left( \lambda \left( D^{2}u_{1}\right) \right) =p_{1}\left( \left\vert
x\right\vert \right) f_{1}\left( u_{2}\right) \text{ for }x\in \mathbb{R}^{N}%
\text{ (}N\geq 3\text{), } \\ 
S_{k}\left( \lambda \left( D^{2}u_{2}\right) \right) =p_{2}\left( \left\vert
x\right\vert \right) f_{2}\left( u_{1}\right) \text{ for }x\in \mathbb{R}^{N}%
\text{ (}N\geq 3\text{).}%
\end{array}%
\right.
\end{equation*}%
To simplify the presentation we wish to mention that our purpose is to
complete and improve all the results in \cite{ZZ} for the more general
system (\ref{11}). By analogy with the work of Zhang and Zhou \cite{ZZ} we
introduce the following notations 
\begin{eqnarray*}
C_{0} &=&(N-1)!/\left[ k_{1}!(N-k_{1})!\right] ,C_{00}=(N-1)!/\left[
k_{2}!(N-k_{2})!\right] , \\
B_{1}^{-}\left( \xi \right) &=&\frac{\xi ^{k_{1}-N}}{C_{0}}e^{-\int_{0}^{\xi
}\frac{1}{C_{0}}t^{k_{1}-1}b_{1}\left( t\right) dt},\text{ }B_{1}^{+}\left(
\xi \right) =\xi ^{N-1}e^{\int_{0}^{\xi }\frac{1}{C_{0}}t^{k_{1}-1}b_{1}%
\left( t\right) dt}p_{1}\left( \xi \right) , \\
B_{2}^{-}\left( \xi \right) &=&\frac{\xi ^{k_{2}-N}}{C_{00}}%
e^{-\int_{0}^{\xi }\frac{1}{C_{00}}t^{k_{2}-1}b_{2}\left( t\right) dt},\text{
}B_{2}^{+}\left( \xi \right) =\xi ^{N-1}e^{\int_{0}^{\xi }\frac{1}{C_{00}}%
t^{k_{2}-1}b_{2}\left( t\right) dt}p_{2}\left( \xi \right) , \\
P_{1}\left( r\right) &=&\int_{0}^{r}\left( B_{1}^{-}\left( r\right)
\int_{0}^{r}B_{1}^{+}\left( t\right) dt\right) ^{\frac{1}{k_{1}}}dr, \\
P_{2}\left( r\right) &=&\int_{0}^{r}\left( B_{2}^{-}\left( r\right)
\int_{0}^{r}B_{2}^{+}\left( t\right) dt\right) ^{\frac{1}{k_{2}}}dr, \\
F_{1,2}\left( r\right) &=&\int_{a_{1}+a_{2}}^{r}\frac{1}{f_{1}^{1/k_{1}}%
\left( t,t\right) +f_{2}^{1/k_{2}}\left( t,t\right) }dt\text{ for }r\geq
a_{1}+a_{2}>0\text{, }a_{1}\geq 0\text{, }a_{1}\geq 0\text{,} \\
P_{1}\left( \infty \right) &=&\lim_{r\rightarrow \infty }P_{1}\left(
r\right) ,\text{ }P_{2}\left( \infty \right) =\lim_{r\rightarrow \infty
}P_{2}\left( r\right) \text{, }F_{1,2}\left( \infty \right)
=\lim_{s\rightarrow \infty }F_{1,2}\left( s\right) .
\end{eqnarray*}%
We will always assume that the variable weights functions $%
b_{1},b_{2},p_{1},p_{2}$ and the nonlinearities $f_{1}$, $f_{2}$ satisfy:

(P1)\quad $b_{1},b_{2}:\left[ 0,\infty \right) \rightarrow \left[ 0,\infty
\right) $ and $p_{1},p_{2}:\left[ 0,\infty \right) \rightarrow \left[
0,\infty \right) $ are spherically symmetric continuous functions (i.e.,%
\textit{\ }$p_{i}\left( x\right) =p_{i}\left( \left\vert x\right\vert
\right) $ and $b_{i}\left( x\right) =b_{i}\left( \left\vert x\right\vert
\right) $ for $i=1,2$)\textit{;}

(C1)\quad $f_{1}$, $f_{2}:\left[ 0,\infty \right) \times \left[ 0,\infty
\right) \rightarrow \left[ 0,\infty \right) $ are continuous and increasing.

Here is a first (somewhat suprizing) result:

\begin{theorem}
\label{th1}We \textit{assume that }$F_{1,2}\left( \infty \right) =\infty $
and \textrm{(P1),} hold\textit{. Furthermore, if }$f_{1}$ and $f_{2}$
satisfy \textrm{(C1),} then the system (\ref{11}) has at least one positive
radial solution $\left( u_{1},u_{2}\right) \in C^{2}\left( \left[ 0,\infty
\right) \right) \times C^{2}\left( \left[ 0,\infty \right) \right) $ with
central value in $\left( a_{1},a_{2}\right) $. Moreover, the following hold:

1.)\quad If $P_{1}\left( \infty \right) +P_{2}\left( \infty \right) <\infty $
then $\lim_{r\rightarrow \infty }u_{1}\left( r\right) <\infty $ and $%
\lim_{r\rightarrow \infty }u_{2}\left( r\right) <\infty .$

2.)\quad If $P_{1}\left( \infty \right) =\infty $ and $P_{2}\left( \infty
\right) =\infty $ then $\lim_{r\rightarrow \infty }u_{1}\left( r\right)
=\infty $ and $\lim_{r\rightarrow \infty }u_{2}\left( r\right) =\infty $.
\end{theorem}

In the same spirit we also have, our next result:

\begin{theorem}
\label{th2}Assume that the hypotheses \textrm{(P1) }and \textrm{(C1)} are
satisfied. \ Then, if $F_{1,2}\left( \infty \right) <\infty $ , $P_{1}\left(
\infty \right) +P_{2}\left( \infty \right) <\infty $ and $P_{1}\left( \infty
\right) +P_{2}\left( \infty \right) <F_{1,2}\left( \infty \right) $ the
system (\ref{11}) has one positive bounded radial solution $\left(
u_{1},u_{2}\right) \in C^{2}\left( \left[ 0,\infty \right) \right) \times
C^{2}\left( \left[ 0,\infty \right) \right) $, with central value in $\left(
a_{1},a_{2}\right) ,$ such that%
\begin{equation*}
\left\{ 
\begin{array}{c}
a_{1}+f_{1}^{1/k_{1}}\left( a_{1},a_{2}\right) P_{1}\left( r\right) \leq
u_{1}\left( r\right) \leq F_{1,2}^{-1}\left( P_{1}\left( r\right)
+P_{2}\left( r\right) \right) , \\ 
a_{2}+f_{2}^{1/k_{2}}\left( a_{1},a_{2}\right) P_{2}\left( r\right) \leq
u_{2}\left( r\right) \leq F_{1,2}^{-1}\left( P_{1}\left( r\right)
+P_{2}\left( r\right) \right) .%
\end{array}%
\right. \text{ }
\end{equation*}
\end{theorem}

\section{Proofs of the main results}

We give in this section the proof of Theorem \ref{th1} and Theorem \ref{th2}%
. For the readers' convenience, we recall the radial form of the $k$-Hessian
operator.

\begin{remark}
If $u:\mathbb{R}^{N}\rightarrow \mathbb{R}$ is radially symmetric then a
calculation show%
\begin{equation*}
S_{k}\left( \lambda \left( D^{2}u\left( r\right) \right) \right)
=r^{1-N}C_{N-1}^{k-1}\left[ \frac{r^{N-k}}{k}\left( u^{^{\prime }}\left(
r\right) \right) ^{k}\right] ^{\prime }\text{ ,}
\end{equation*}%
where the prime denotes differentiation with respect to $r=\left\vert
x\right\vert $ and $C_{N-1}^{k-1}=(N-1)!/\left[ (k-1)!(N-k)!\right] $.
\end{remark}

This Remark is well known, see for example \cite{BAO} or \cite{SA}.

\subparagraph{\textbf{Proof of the Theorems \protect\ref{th1} and \protect
\ref{th2}. }}

We start by showing that the system (\ref{11}) has positive radial
solutions. For this purpose, we show that the ordinary differential
equations system 
\begin{equation}
\left\{ 
\begin{array}{l}
C_{N-1}^{k_{1}-1}\left[ \frac{r^{N-k_{1}}}{k_{1}}e^{\int_{0}^{r}\frac{1}{%
C_{0}}t^{k_{1}-1}b_{1}\left( t\right) dt}\left( u_{1}^{^{\prime }}\left(
r\right) \right) ^{k_{1}}\right] ^{\prime }=r^{N-1}e^{\int_{0}^{r}\frac{1}{%
C_{0}}t^{k_{1}-1}b_{1}\left( t\right) dt}p_{1}\left( r\right) f_{1}\left(
u_{1}\left( r\right) ,u_{2}\left( r\right) \right) \text{ for }r>0, \\ 
C_{N-1}^{k_{2}-1}\left[ \frac{r^{N-k_{2}}}{k_{2}}e^{\int_{0}^{r}\frac{1}{%
C_{00}}t^{k_{2}-1}b_{2}\left( t\right) dt}\left( u_{2}^{^{\prime }}\left(
r\right) \right) ^{k_{2}}\right] ^{\prime }=r^{N-1}e^{\int_{0}^{r}\frac{1}{%
C_{00}}t^{k_{2}-1}b_{2}\left( t\right) dt}p_{2}\left( r\right) f_{2}\left(
u_{1}\left( r\right) ,u_{2}\left( r\right) \right) \text{ for }r>0, \\ 
u_{1}^{\prime }\left( r\right) \geq 0\text{ and }u_{2}^{\prime }\left(
r\right) \geq 0\text{ for }r\in \left[ 0,\infty \right) , \\ 
u_{1}\left( 0\right) =a_{1}\text{ and }u_{2}\left( 0\right) =a_{2},%
\end{array}%
\right.  \label{ss}
\end{equation}%
has solution. Therefore, at least one solution of (\ref{ss}) can be obtained
using successive approximation by defining the sequences $\left\{
u_{1}^{m}\right\} ^{m\geq 1}$ and $\left\{ u_{2}^{m}\right\} ^{m\geq 1}$ on $%
\left[ 0,\infty \right) $ in the following way $u_{1}^{0}=a_{1},$ $%
u_{2}^{0}=a_{2}$ for $r\geq 0$ and 
\begin{equation}
\left\{ 
\begin{array}{c}
u_{1}^{m}\left( s\right) =a_{1}+\int_{0}^{r}\left[ B_{1}^{-}\left( t\right)
\int_{0}^{t}B_{1}^{+}\left( s\right) f_{1}\left( u_{1}^{m-1}\left( s\right)
,u_{2}^{m-1}\left( s\right) \right) ds\right] ^{1/k_{1}}dt, \\ 
u_{2}^{m}\left( s\right) =a_{2}+\int_{0}^{r}\left[ B_{2}^{-}\left( t\right)
\int_{0}^{t}B_{2}^{+}\left( s\right) f_{2}\left( u_{1}^{m-1}\left( s\right)
,u_{2}^{m-1}\left( s\right) \right) ds\right] ^{1/k_{2}}dt.%
\end{array}%
\right.  \label{recs}
\end{equation}%
It is easy to see that $\left\{ u_{1}^{m}\right\} ^{m\geq 1}$ and $\left\{
u_{2}^{m}\right\} ^{m\geq 1}$ are non-decreasing on $\left[ 0,\infty \right) 
$. Indeed, we consider 
\begin{eqnarray*}
u_{1}^{1}\left( r\right) &=&a_{1}+\int_{0}^{r}\left[ B_{1}^{-}\left(
t\right) \int_{0}^{t}B_{1}^{+}\left( s\right) f_{1}\left( u_{1}^{0}\left(
s\right) ,u_{2}^{0}\left( s\right) \right) ds\right] ^{1/k_{1}}dt \\
&=&a_{1}+\int_{0}^{r}\left[ B_{1}^{-}\left( t\right)
\int_{0}^{t}B_{1}^{+}\left( s\right) f_{1}\left( a_{1},a_{2}\right) ds\right]
^{1/k_{1}}dt \\
&\leq &a_{1}+\int_{0}^{r}\left[ B_{1}^{-}\left( t\right)
\int_{0}^{t}B_{1}^{+}\left( s\right) f_{1}\left( u_{1}^{1}\left( s\right)
,u_{2}^{1}\left( s\right) \right) ds\right] ^{1/k_{1}}dt=u_{1}^{2}\left(
r\right) .
\end{eqnarray*}%
This implies that 
\begin{equation*}
u_{1}^{1}\left( r\right) \leq u_{1}^{2}\left( r\right) \text{ which further
produces }u_{1}^{2}\left( r\right) \leq u_{1}^{3}\left( r\right) .
\end{equation*}%
Continuing, an induction argument applied to (\ref{recs}) show that for any $%
r\geq 0$ we have 
\begin{equation*}
u_{1}^{m}\left( r\right) \leq u_{1}^{m+1}\left( r\right) \text{ and }%
u_{2}^{m}\left( r\right) \leq u_{2}^{m+1}\left( r\right) \text{ for any }%
m\in \mathbb{N}
\end{equation*}%
i.e.,$\left\{ u_{1}^{m}\right\} ^{m\geq 1}$ and $\left\{ u_{2}^{m}\right\}
_{j=1,2}^{m\geq 1}$ are non-decreasing on $\left[ 0,\infty \right) $. By the
monotonicity of $\left\{ u_{1}^{m}\right\} ^{m\geq 1}$ and $\left\{
u_{2}^{m}\right\} ^{m\geq 1}$ we have the inequalities%
\begin{eqnarray}
C_{N-1}^{k_{1}-1}\left\{ \frac{r^{N-k_{1}}}{k_{1}}e^{\int_{0}^{r}\frac{1}{%
C_{0}}t^{k_{1}-1}b_{1}\left( t\right) dt}\left[ \left( u_{1}^{m}\left(
r\right) \right) ^{\prime }\right] ^{k_{1}}\right\} ^{\prime } &\leq
&B_{1}^{+}\left( r\right) f_{1}\left( u_{1}^{m}\left( r\right)
,u_{2}^{m}\left( r\right) \right) ,  \label{gen1} \\
C_{N-1}^{k_{2}-1}\left\{ \frac{r^{N-k_{2}}}{k_{2}}e^{\int_{0}^{r}\frac{1}{%
C_{00}}t^{k_{2}-1}b_{2}\left( t\right) dt}\left[ \left( u_{2}^{^{m}}\left(
r\right) \right) ^{\prime }\right] ^{k_{2}}\right\} ^{\prime } &\leq
&B_{2}^{+}\left( r\right) f_{2}\left( u_{1}^{m}\left( r\right)
,u_{2}^{m}\left( r\right) \right) .  \label{gen2}
\end{eqnarray}%
After integration from $0$ to $r$, an easy calculation yields 
\begin{eqnarray}
\left( u_{1}^{m}\left( r\right) \right) ^{\prime } &\leq &\left(
B_{1}^{-}\left( r\right) \int_{0}^{r}B_{1}^{+}\left( t\right) f_{1}\left(
u_{1}^{m}\left( t\right) ,u_{2}^{m}\left( t\right) \right) dt\right) ^{\frac{%
1}{k_{1}}}  \notag \\
&\leq &\left( B_{1}^{-}\left( r\right) \int_{0}^{r}B_{1}^{+}\left( t\right)
f_{1}\left( u_{1}^{m}\left( t\right) +u_{2}^{m}\left( t\right)
,u_{1}^{m}\left( t\right) +u_{2}^{m}\left( t\right) \right) dt\right) ^{%
\frac{1}{k_{1}}}  \label{exin} \\
&\leq &\left( f_{1}^{1/k_{1}}+f_{2}^{1/k_{2}}\right) \left( u_{1}^{m}\left(
r\right) +u_{2}^{m}\left( r\right) ,u_{1}^{m}\left( r\right)
+u_{2}^{m}\left( r\right) \right) \left( B_{1}^{-}\left( r\right)
\int_{0}^{r}B_{1}^{+}\left( t\right) dt\right) ^{\frac{1}{k_{1}}}.  \notag
\end{eqnarray}%
As before, exactly the same type of conclusion holds for $\left(
u_{2}^{m}\left( r\right) \right) ^{\prime }$: 
\begin{eqnarray}
\left( u_{2}^{m}\left( r\right) \right) ^{\prime } &\leq &\left(
B_{2}^{-}\left( r\right) \int_{0}^{r}B_{2}^{+}\left( z\right) f_{2}\left(
u_{1}^{m}\left( z\right) ,u_{2}^{m}\left( z\right) \right) dz\right)
^{1/k_{2}}  \notag \\
&\leq &\left( f_{1}^{1/k_{1}}+f_{2}^{1/k_{2}}\right) \left( u_{1}^{m}\left(
r\right) +u_{2}^{m}\left( r\right) ,u_{1}^{m}\left( r\right)
+u_{2}^{m}\left( r\right) \right) \left( B_{2}^{-}\left( r\right)
\int_{0}^{r}B_{2}^{+}\left( t\right) dt\right) ^{\frac{1}{k_{2}}}.
\label{exin2}
\end{eqnarray}%
Summing the two previous inequalities (\ref{exin}) and (\ref{exin2}), we
obtain 
\begin{equation}
\frac{\left( u_{1}^{m}\left( r\right) +u_{2}^{m}\left( r\right) \right)
^{\prime }}{\left( f_{1}^{1/k_{1}}+f_{2}^{1/k_{2}}\right) \left(
u_{1}^{m}\left( r\right) +u_{2}^{m}\left( r\right) ,u_{1}^{m}\left( r\right)
+u_{2}^{m}\left( r\right) \right) }\leq \left( B_{1}^{-}\left( r\right)
\int_{0}^{r}B_{1}^{+}\left( t\right) dt\right) ^{\frac{1}{k_{1}}}+\left(
B_{2}^{-}\left( r\right) \int_{0}^{r}B_{2}^{+}\left( t\right) dt\right) ^{%
\frac{1}{k_{2}}}.  \label{mat}
\end{equation}%
Integrating from $0$ to $r$ the inequality (\ref{mat}), we obtain%
\begin{equation*}
\int_{a_{1}+a_{2}}^{u_{1}^{m}\left( r\right) +u_{2}^{m}\left( r\right) }%
\frac{1}{f_{1}^{1/k_{1}}\left( t,t\right) +f_{2}^{1/k_{2}}\left( t,t\right) }%
dt\leq P_{1}\left( r\right) +P_{2}\left( r\right) .
\end{equation*}%
We now have 
\begin{equation}
F_{1,2}\left( u_{1}^{m}\left( r\right) +u_{2}^{m}\left( r\right) \right)
\leq P_{1}\left( r\right) +P_{2}\left( r\right) ,  \label{ints}
\end{equation}%
which will play a basic role in the proof of our main results. The
inequalities (\ref{ints}) can be reformulated as 
\begin{equation}
u_{1}^{m}\left( r\right) +u_{2}^{m}\left( r\right) \leq F_{1,2}^{-1}\left(
P_{1}\left( r\right) +P_{2}\left( r\right) \right) .  \label{int}
\end{equation}%
This can be easily seen from the fact that $F_{1,2}$ is a bijection with the
inverse function $F_{1,2}^{-1}$ strictly increasing on $\left[ 0,\infty
\right) $. So, we have found upper bounds for $\left\{ u_{1}^{m}\right\}
^{m\geq 1}$ and $\left\{ u_{2}^{m}\right\} ^{m\geq 1}$ which are dependent
of $r$. We are now ready to give a complete proof of the Theorems \ref{th1}-%
\ref{th2}.

\textbf{Proof of Theorem \ref{th1} completed:} When $F_{1,2}\left( \infty
\right) =\infty $ it follows that the sequences $\left\{ u_{j}^{m}\right\}
_{j=1,2}^{m\geq 1}$ are bounded and equicontinuous on $\left[ 0,c_{0}\right] 
$ for arbitrary $c_{0}>0$. Possibly after passing to a subsequence, we may
assume that$\left\{ u^{m}\right\} _{j=1,2}^{m\geq 1}$ converges uniformly to 
$\left\{ u_{j}\right\} _{j=1,2}$ on $\left[ 0,c_{0}\right] $. At the end of
this process, we conclude by the arbitrariness of $c_{0}>0$, that $\left(
u_{1},u_{2}\right) $ is positive entire solution of system (\ref{ss}). The
solution constructed in this way will be radially symmetric. Since the
radial solutions of the ordinary differential equations system (\ref{ss})
are solutions (\ref{11}) it follows that the radial solutions of (\ref{11})
with $u_{1}\left( 0\right) =a_{1},$ $u_{2}\left( 0\right) =a_{2}$ satisfy: 
\begin{eqnarray}
u_{1}\left( r\right) &=&a_{1}+\int_{0}^{r}\left( B_{1}^{-}\left( y\right)
\int_{0}^{y}B_{1}^{+}\left( t\right) f_{1}\left( u_{1}\left( t\right)
,u_{2}\left( t\right) \right) dt\right) ^{1/k_{1}}dy,\text{ }r\geq 0,
\label{eq1} \\
u_{2}\left( r\right) &=&a_{2}+\int_{0}^{r}\left( B_{2}^{-}\left( y\right)
\int_{0}^{y}B_{2}^{+}\left( t\right) f_{2}\left( u_{1}\left( t\right)
,u_{2}\left( t\right) \right) dt\right) ^{1/k_{2}}dy,\text{ }r\geq 0.
\label{eq2}
\end{eqnarray}%
Next, it is easy to verify that the Cases 1. and 2. occur.

\textbf{Case 1.):} When $P_{1}\left( \infty \right) +P_{2}\left( \infty
\right) <\infty ,$ it is not difficult to deduce from (\ref{eq1}) and (\ref%
{eq2}) that%
\begin{equation*}
u_{1}\left( r\right) +u_{2}\left( r\right) \leq F_{1,2}^{-1}\left(
P_{1}\left( \infty \right) +P_{2}\left( \infty \right) \right) <\infty \text{
for all }r\geq 0,
\end{equation*}%
and so $\left( u_{1},u_{2}\right) $ is bounded. We next consider:

\textbf{Case 2.):} In the case $P_{1}\left( \infty \right) =P_{2}\left(
\infty \right) =\infty $, we observe that 
\begin{eqnarray}
u_{1}\left( r\right) &=&a_{1}+\int_{0}^{r}\left( B_{1}^{-}\left( t\right)
\int_{0}^{t}B_{1}^{+}\left( s\right) f_{1}\left( u_{1}\left( s\right)
,u_{2}\left( s\right) \right) ds\right) ^{\frac{1}{k_{1}}}dt  \notag \\
&\geq &a_{1}+f_{1}^{1/k_{1}}\left( a_{1},a_{2}\right) \int_{0}^{r}\left(
B_{1}^{-}\left( t\right) \int_{0}^{t}B_{1}^{+}\left( s\right) ds\right) ^{%
\frac{1}{k_{1}}}dt  \label{i1} \\
&=&a_{1}+f_{1}^{1/k_{1}}\left( a_{1},a_{2}\right) P_{1}\left( r\right) . 
\notag
\end{eqnarray}%
The same computations as in (\ref{i1}) yields%
\begin{equation*}
u_{2}\left( r\right) \geq a_{2}+f_{2}^{1/k_{2}}\left( a_{1},a_{2}\right)
P_{2}\left( r\right) .
\end{equation*}%
and passing to the limit as $r\rightarrow \infty $ in (\ref{i1}) and in the
above inequality we conclude that%
\begin{equation*}
\lim_{r\rightarrow \infty }u_{1}\left( r\right) =\lim_{r\rightarrow \infty
}u_{2}\left( r\right) =\infty ,
\end{equation*}%
which yields the result. We now turn to:

\textbf{Proof of Theorem \ref{th2} completed: }

In view of the above analysis, the proof can be easily deduced from 
\begin{equation*}
F_{1,2}\left( u_{1}^{m}\left( r\right) +u_{2}^{m}\left( r\right) \right)
\leq P_{1}\left( \infty \right) +P_{2}\left( \infty \right) <F_{1,2}\left(
\infty \right) <\infty ,
\end{equation*}%
Indeed, since $F_{1,2}^{-1}$ is strictly increasing on $\left[ 0,\infty
\right) $, we find that%
\begin{equation*}
u_{1}^{m}\left( r\right) +u_{2}^{m}\left( r\right) \leq F_{1,2}^{-1}\left(
P_{1}\left( \infty \right) +P_{2}\left( \infty \right) \right) <\infty ,
\end{equation*}%
and then the non-decreasing sequences $\left\{ u_{1}^{m}\left( r\right)
\right\} ^{m\geq 1}$ and $\left\{ u_{2}^{m}\left( r\right) \right\} ^{m\geq
1}$ are bounded above for all $r\geq 0$ and all $m$. The final step, is to
conclude that $\left( u_{1}^{m}\left( r\right) ,u_{2}^{m}\left( r\right)
\right) \rightarrow \left( u_{1}\left( r\right) ,u_{2}\left( r\right)
\right) $ as $m\rightarrow \infty $ and the limit functions $u_{1}$ and $%
u_{2}$ are positive entire bounded radial solutions of system (\ref{11}).%
\textbf{\ }

\begin{remark}
Make the same assumptions as in Theorem \ref{th1} or Theorem \ref{th2} on $%
b_{1}$, $p_{1}$, $f_{1}$, $b_{2}$, $p_{2}$, $f_{2}$. If, in addition,%
\begin{eqnarray}
p_{1}\left( \left\vert x\right\vert \right) &\geq &\left( C_{N-1}^{k_{1}-1}%
\frac{N-k_{1}}{k_{1}}\left\vert x\right\vert ^{-N}-b_{1}\left( \left\vert
x\right\vert \right) \left\vert x\right\vert ^{k_{1}-N}\right)
\int_{0}^{\left\vert x\right\vert }\frac{s^{N-1}}{C_{0}}p_{1}\left( s\right)
ds\text{ for every }x\in \mathbb{R}^{N},  \label{ing1} \\
p_{2}\left( \left\vert x\right\vert \right) &\geq &\left( C_{N-1}^{k_{2}-1}%
\frac{N-k_{2}}{k_{2}}\left\vert x\right\vert ^{-N}-b_{2}\left( \left\vert
x\right\vert \right) \left\vert x\right\vert ^{k_{2}-N}\right)
\int_{0}^{\left\vert x\right\vert }\frac{s^{N-1}}{C_{00}}p_{2}\left(
s\right) ds\text{ for every }x\in \mathbb{R}^{N},  \label{ing2}
\end{eqnarray}%
then, the solution $\left( u_{1},u_{2}\right) $ is convex.
\end{remark}

\subparagraph{Proof.}

It is clear that%
\begin{equation}
C_{N-1}^{k_{1}-1}\left[ \frac{r^{N-k_{1}}}{k_{1}}e^{\int_{0}^{r}\frac{1}{%
C_{0}}t^{k_{1}-1}b_{1}\left( t\right) dt}\left( u_{1}^{^{\prime }}\left(
r\right) \right) ^{k_{1}}\right] ^{\prime }=r^{N-1}e^{\int_{0}^{r}\frac{1}{%
C_{0}}t^{k_{1}-1}b_{1}\left( t\right) dt}p_{1}\left( r\right) f_{1}\left(
u_{1}\left( r\right) ,u_{2}\left( r\right) \right) ,  \label{jo}
\end{equation}%
and integrating from $0$ to $r$ yields%
\begin{eqnarray*}
r^{N-k_{1}}e^{\int_{0}^{r}\frac{1}{C_{0}}t^{k_{1}-1}b_{1}\left( t\right)
dt}\left( u_{1}^{^{\prime }}\left( r\right) \right) ^{k_{1}} &=&\int_{0}^{r}%
\frac{s^{N-1}}{C_{0}}e^{\int_{0}^{s}\frac{1}{C_{0}}t^{k_{1}-1}b_{1}\left(
t\right) dt}p_{1}\left( s\right) f_{1}\left( u_{1}\left( s\right)
,u_{2}\left( s\right) \right) ds \\
&\leq &f_{1}\left( u_{1}\left( r\right) ,u_{2}\left( r\right) \right)
\int_{0}^{r}\frac{s^{N-1}}{C_{0}}e^{\int_{0}^{s}\frac{1}{C_{0}}%
t^{k_{1}-1}b_{1}\left( t\right) dt}p_{1}\left( s\right) ds,
\end{eqnarray*}%
which produces%
\begin{eqnarray}
\left( \frac{u_{1}^{^{\prime }}\left( r\right) }{r}\right) ^{k_{1}} &\leq
&f_{1}\left( u_{1}\left( r\right) ,u_{2}\left( r\right) \right)
r^{-N}e^{-\int_{0}^{r}\frac{1}{C_{0}}t^{k_{1}-1}b_{1}\left( t\right)
dt}\int_{0}^{r}\frac{s^{N-1}}{C_{0}}e^{\int_{0}^{s}\frac{1}{C_{0}}%
t^{k_{1}-1}b_{1}\left( t\right) dt}p_{1}\left( s\right) ds  \notag \\
&\leq &f_{1}\left( u_{1}\left( r\right) ,u_{2}\left( r\right) \right)
r^{-N}\int_{0}^{r}\frac{s^{N-1}}{C_{0}}p_{1}\left( s\right) ds.  \label{conv}
\end{eqnarray}%
On the other hand the inequality (\ref{jo}) can be writen in the following
way%
\begin{equation}
C_{N-1}^{k_{1}-1}u_{1}^{^{\prime \prime }}\left( r\right) \left( \frac{%
u_{1}^{^{\prime }}\left( r\right) }{r}\right) ^{k_{1}-1}+C_{N-1}^{k_{1}-1}%
\frac{N-k_{1}}{k_{1}}\left( \frac{u_{1}^{^{\prime }}\left( r\right) }{r}%
\right) ^{k_{1}}+b_{1}\left( r\right) \left( u_{1}^{^{\prime }}\left(
r\right) \right) ^{k_{1}}=p_{1}\left( r\right) f_{1}\left( u_{1}\left(
r\right) ,u_{2}\left( r\right) \right) .  \label{fc}
\end{equation}%
Using the inequality (\ref{conv}) in (\ref{fc}) we obtain%
\begin{eqnarray*}
p_{1}\left( r\right) f_{1}\left( u_{1}\left( r\right) ,u_{2}\left( r\right)
\right) \leq C_{N-1}^{k_{1}-1}u_{1}^{^{\prime \prime }}\left( r\right)
\left( \frac{u_{1}^{^{\prime }}\left( r\right) }{r}\right) ^{k_{1}-1}
&&+C_{N-1}^{k_{1}-1}\frac{N-k_{1}}{k_{1}}f_{1}\left( u_{1}\left( r\right)
,u_{2}\left( r\right) \right) r^{-N}\int_{0}^{r}\frac{s^{N-1}}{C_{0}}%
p_{1}\left( s\right) ds \\
&&+b_{1}\left( r\right) r^{k_{1}}f_{1}\left( u_{1}\left( r\right)
,u_{2}\left( r\right) \right) r^{-N}\int_{0}^{r}\frac{s^{N-1}}{C_{0}}%
p_{1}\left( s\right) ds,
\end{eqnarray*}%
from which we have%
\begin{equation*}
f_{1}\left( u_{1}\left( r\right) ,u_{2}\left( r\right) \right) \left[
p_{1}\left( r\right) -\left( C_{N-1}^{k_{1}-1}\frac{N-k_{1}}{k_{1}}%
r^{-N}-b_{1}\left( r\right) r^{k_{1}-N}\right) \int_{0}^{r}\frac{s^{N-1}}{%
C_{0}}p_{1}\left( s\right) ds\right] \leq C_{N-1}^{k_{1}-1}u_{1}^{^{\prime
\prime }}\left( r\right) \left( \frac{u_{1}^{^{\prime }}\left( r\right) }{r}%
\right) ^{k_{1}-1},
\end{equation*}%
which completes the proof of $u_{1}^{^{\prime \prime }}\left( r\right) \geq
0 $. A similar argument produces $u_{2}^{^{\prime \prime }}\left( r\right)
\geq 0$.\ We also remark that, in the simple case $b_{1}=b_{2}=0,$ $%
s^{N-1}p_{1}\left( s\right) $ and $s^{N-1}p_{2}\left( s\right) $ are
increasing then (\ref{ing1}) and (\ref{ing2}) hold.

\textbf{Acknowledgement. }The author would like to thank to the editors and
reviewers for valuable comments and suggestions which contributed to improve
this article.

\end{document}